\documentclass[pre,twocolumn,superscriptaddress]{revtex4}  

\usepackage[paperwidth=210mm,paperheight=297mm,centering,hmargin=2cm,vmargin=2.5cm]{geometry}

\usepackage[pdftex]{graphicx}
\usepackage{amsmath,amsfonts,amssymb}
\usepackage{morefloats}
\usepackage{xcolor} 
\usepackage{tabularx}

\usepackage{newfloat,algcompatible}
\usepackage{etoolbox}
\DeclareFloatingEnvironment[
    fileext=loa,
    listname=List of Algorithms,
    name= Algorithm,
    placement=tbhp,
]{algorithm}

\makeatletter
\renewcommand{\fnum@algorithm}{\small\textbf{\algorithmname~\thealgorithm}}
\makeatother

\definecolor{gris25}{gray}{0.5}
\definecolor{gris5}{gray}{0.7}
\definecolor{gris75}{gray}{0.9}
\definecolor{green2}{rgb}{0,0.7,0.3}
\definecolor{bluemoi}{rgb}{0.25,0.50 ,0.75} 

\usepackage{hyperref}
\hypersetup{
  bookmarksopen=false,
  pdftitle=" ",
  pdfauthor=" ", 
  pdfsubject=" ", 
  pdftoolbar=true, 
  pdfmenubar=true, 
  pdfhighlight=/O, 
  colorlinks=true, 
  pdfpagemode=UseNone, 
  pdfpagelayout=SinglePage, 
  pdffitwindow=true, 
  linkcolor=bluemoi, 
  citecolor=bluemoi, 
  urlcolor=bluemoi 
}

\makeatletter
\renewcommand{\figurename}{Figure}
\renewcommand{\fnum@figure}{\small\textbf{\figurename~\thefigure}}
\renewcommand{\thefigure}{\arabic{figure}}
\makeatother

\makeatletter
\renewcommand{\tablename}{Table}
\renewcommand{\fnum@table}{\small\textbf{\tablename~\thetable}}
\renewcommand{\thetable}{\arabic{table}}
\makeatother

\begin{document}

\title{Generating French virtual commuting networks at municipality level}

\author{Maxime Lenormand}\affiliation{IRSTEA, LISC, 24 avenue des Landais, 63172 AUBIERE, France}
\author{Sylvie Huet}\affiliation{IRSTEA, LISC, 24 avenue des Landais, 63172 AUBIERE, France}
\author{Floriana Gargiulo}\affiliation{IRSTEA, LISC, 24 avenue des Landais, 63172 AUBIERE, France}

\begin{abstract} 
We aim to generate virtual commuting networks in the rural regions of France in order to study the dynamics of their municipalities. Since it will be necessary to model small commuting flows between municipalities with a few hundred or thousand inhabitants, we have opted for the stochastic model presented by \cite{Gargiulo2012}. This model reproduces various possible complete networks using an iterative process, stochastically selecting a workplace in the region for each commuter living in the municipality of a region. The choice is made considering the job offers in each municipality of the region and the distance to all of the possible destinations. This paper will present methods for adapting and implementing this model to generate commuting networks between municipalities for regions in France. We address three different issues: How can we generate a reliable virtual commuting network for a region that is highly dependent on other regions for the satisfaction of its resident's demands for employment? What about a convenient deterrence function? How to calibrate the model when detailed data is not available? Our solution proposes an extended job search geographical base for commuters living in the municipalities, we compare two different deterrence functions and we show that the parameter is a constant for network linking municipalities in France. 
\end{abstract}

\maketitle

\section{Introduction}

The connection between the home and workplace plays a central role in understanding the socio-economic relations in a network of rural municipalities \cite{Clark2003,Reggiani2010}. Indeed, new economic theories assume local positive dynamics can be explained by implicit geographical money transfers made by commuters or retired people (see for example \cite{Davezies2009}). Simulation is becoming an increasingly convenient tool to study populations and their interactions over the space. That is particularly the case with the individual-based approaches which allow studying theories at the individual level since they simulate the variations in how individuals interact with each other and with their environment. Recent modeling reviews show the increasing use of such a tool \cite{Parker2003,Waddell2003,Bousquet2004,Verburg2004,Rindfuss2004,Birkin2012}. However, these approaches require generation models capable of building reliable virtual commuting networks that consider each individual within a population. That is the case in the SimVillages dynamic micro-simulation model we developed during the PRIMA project\footnote{PRototypical policy Impacts on Multifunctional Activities in rural municipalities - EU 7th Framework Research Programme; 2008-2011; \url{https://prima.cemagref.fr/the-project}}. Indeed, in the SimVillages model, after generating a synthetic population of individuals \cite{Gargiulo2010}, it is necessary to choose a place of work for each worker within this population because a commuting origin-destination table was unavailable.

The goal of the European PRIMA project was to understand the dynamics of rural municipalities in France. 95\% of them have less than 3000
inhabitants. This means that most of the commuting flows we want to study are weak, with a spatial distribution very difficult to predict with the available variables at an aggregated level. This is why we opt for the stochastic model recently proposed by \cite{Gargiulo2012}. Moreover, we want to consider the commuting network on different dates. Detailed data regarding flows between pairs of municipalities are only available in France for the year 1999. For other dates, the only reliable data is aggregated data for each municipality, which describes how many people work outside of the municipality and how many come from outside of the municipality to work. Such data lacks precision regarding the various places of work and the various municipalities where citizens reside. Then we also choose the \cite{Gargiulo2012} model for its ability to generate a population of individuals on a commuting network, starting from this data. This model reproduces the complete network using an iterative process that stochastically selects a workplace in the region for each commuter living in the municipality of the region. The choice is made while considering the job offers in each municipality of the region and the distance to all possible destinations. It differs from the classical generation models presented in \cite{Ortuzar2011} since it is a discrete choice model where the individual decision function is inspired by the gravity law model, which is not usually employed on an individual level \cite{Haynes1984,Ortuzar2011,Barthelemy2011}. Moreover, such a model ensures that for every municipality the virtual total numbers of commuters both coming in and going out are the same as the ones supplied by the data.  
This paper presents a method to adapt and implement this model to generate commuting networks between municipalities for regions in France. This implementation has forced us to address three different issues: How can we generate a reliable virtual commuting network for a region highly dependent of other regions to satisfy the need for job for the people living in the municipalities? What about a convenient deterrence function? How should the model be calibrated when detailed data is not available?  

The first problem to solve involves the fact that regions in France are not islands, as presented in the example of \cite{DeMontis2007,DeMontis2010}. Indeed, some of the inhabitants, especially those living close to the borders of the region, are likely to work in municipalities located outside the region of residence. This part, especially if it is significant, causes the generated network to register false if we only consider that people living in the region also work in the region. A method for solving this problem involves generating the commuting network only for people living and working in the region. However, in order to do this it is required that the modeler know the quantity and the place of residence for individuals who work outside but live in the region. Data providing this information is very rare. Therefore, we address this issue by extending the job search geographical base for commuters living in the municipalities to a sufficiently large number of municipalities located outside the region of residence. Then, we compare the model without outside municipalities and the model with outside municipalities in 23 regions in France and come to a conclusion regarding the quality of our solution. 

The second problem relates to the form of the deterrence function which governs the impact of distance on choice of the place of work relative to the quantity of job offers. The initial work done by \cite{Gargiulo2012} propose the use of a power law. However, \cite{Barthelemy2011} states that the form of the deterrence function varies greatly, and can sometimes be inspired by an exponential function, such as in \cite{Balcan2009}, or by a power law function as in \cite{Viboud2006}. To choose the much more convenient deterrence function, we have compared the quality of generated networks for 34 regions in France obtained with both the exponential law and the power law. Better results were obtained with the exponential law.

The final problem was related to calibration. The generation model, as with most of the currently used commuting network generation models, has one parameter to calibrate. This parameter governs the impact of distance on the individual decision regarding the place of work relative to the quantity of job offers. This parameter was calibrated through minimization of the Kolmokorov-Smirnov distance between the observed and simulated commuting distance distribution for individuals of the studied region. When detailed data is not available, it is necessary to find a way to determine this parameter. The only available distance that can be used is the Euclidian distance. While detailed commuting network data was available for the year 1999 and could be used for calibration, it was not available for earlier or more recent years. Though it may be possible to assume the parameter value does not change over time, a transportation network can evolve greatly at the local level to reduce the time distance. Such a change cannot be recorded when using the Euclidian distance. A solution was finally found. Using 34 regions in France, we show that every region can be generated using a constant value for the parameter. Then, we assume that the parameter value is constant over time and space.

\section{Material and methods}

\subsection{The French case-studies and data from the French statistical office \label{observeddata}}

A complete description of the regions from which the network was generated is provided in Table \ref{tabdata}. These regions have been randomly chosen for their diversity in terms of number of municipalities, number of commuters and surface areas. Some correspond to an administrative region of France while others are closer to the county (known as "departements", a French administrative unit). These two types of case studies are called "region" hereafter.

The French Statistical Office ($INSEE$) collects information regarding each individual's residence and place of work. From this collected data, the Maurice Halbwachs Center or the $INSEE$ make the following data available for every researcher:

\begin{enumerate}
\item in 1999, the observed commuting networks i.e. data regarding the numbers of individuals commuting from location $i$ to location $j$ for every municipality of a region (called "observed data" hereafter);
\item in 1999, the total number of commuters, the total job offers and the total number of workers in residence for every municipality. These data allow computations to be made for the number of workers that commute to their office of employment for each municipality. 
\item The Lambert coordinates for each municipality are easy to find on the internet. They are used to compute the Euclidian distance between each pair of municipalities.
\end{enumerate}

We used the data sets 2 and 3 as inputs of the algorithms described in this paper to simulate commuting networks (noted $S$). We compare these simulated commuting networks to "real" network (noted $R$) built from the observed data of the data set 1.

\subsection{The Gargiulo's model}

Consider a region composed of $n$ municipalities. We can model the observed commuting network starting from matrix $R \in \mathrm{M}_{n\times n}(\mathbb{N})$ where $R_{ij}$ represents the number of commuters from municipality $i$ (in the region) to municipality $j$ (in the region). This matrix represents the light gray origin-destination table presented in Table \ref{ODTable1}.

The inputs of the algorithm are: 

\begin{itemize}
\item $D=(d_{ij})_{1 \leq i,j \leq n}$ the Euclidean distance matrix between municipalities.
\item $I_j$ the number of in-commuters from the region to municipality $j$ of the region, $1 \leq j \leq n$ (i.e. the number of individuals living in the region in municipality $i$ ($i\neq j$) and working in municipality $j$).
\item $O_i$ the number of out-commuters from municipality $i$ of the region to the region, $1 \leq i \leq n$ (i.e. the number of individuals working in the region in municipality $j$ ($j\neq i$) and living in municipality $i$).
\end{itemize}

$I_k$ and $O_k$ can be respectively assimilated to the job offers for those employed in the region and the job demand of those employed in the region for municipality $k$, $1 \leq k \leq n$. The algorithm starts with: 

\begin {equation}
I_{j}= \sum_{i=1}^n R_{ij}
\end{equation}

and

\begin {equation}
O_{i}= \sum_{j=1}^n R_{ij}
\end{equation}

The purpose of the model is to generate the light gray origin-destination sub-table of the region described in Table \ref{ODTable1}. To do this it generates matrix $S \in \mathrm{M}_{n\times n}(\mathbb{N})$ where $S_{ij}$ represents the number of commuters from municipality $i$ (in the region) to municipality $j$ (in the region). It's important to note that $S_{ij}=0$ if $i=j$. 
The algorithm assigns to each individual a place of work with a probability based on the distance from the place of residence to every possible place of work and their corresponding job offer. The number of in-commuters for municipality $j$ and the number of out-commuters for municipality $i$ decrease each time an individual living in $i$ is assigned municipality $j$ as a workplace. The algorithm is stopped when all out-commuters have a place of work. The algorithm is described in Algorithm \ref{commut} with $m=n$.

\begin{table*}
	\caption[Origin-destination table for the region]{Origin-destination table for the region; The light gray table represents the commuters living (place of residence RP) and working (place of work WP) in the region for each municipality of the region; The dark gray line represents the number of out-commuters from municipality of the region to the region for each municipality of the region (i.e. the row totals of the light gray table); The dark gray column represents the number of in-commuters from the region to a municipality of the region for each municipality of the region (i.e. the column totals of the light gray table).}
	\label{ODTable1}
		\begin{center}
    \begin{tabular}{c}
			\includegraphics[scale=0.35]{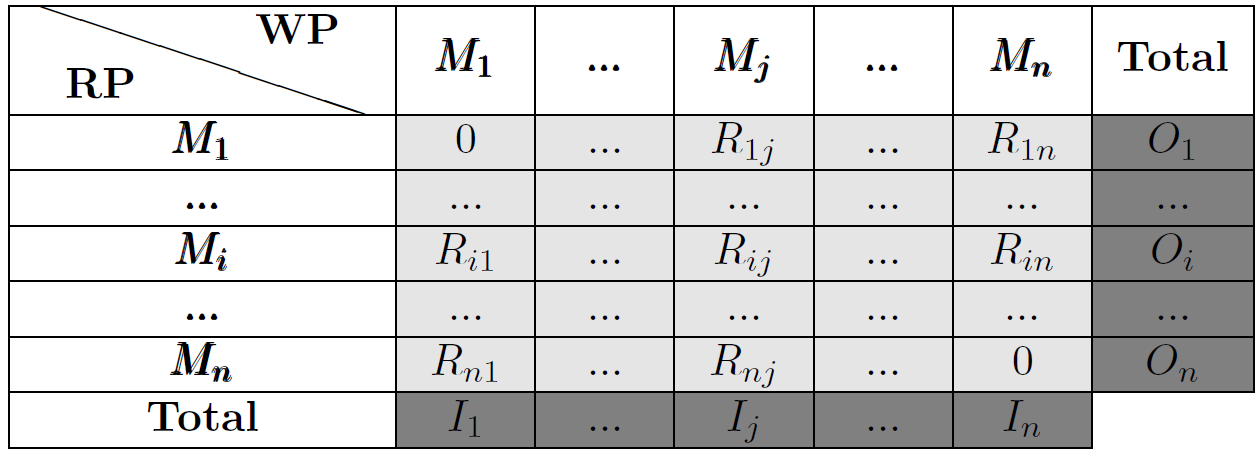}
		\end{tabular}
	\end{center}
\end{table*}

\begin{algorithm}
  {\hrulefill}
	\vspace*{-0.3cm}
	\caption{Commuting generation model}
	{\vspace*{-0.25cm}\hrulefill}
	\label{commut}
	\begin{algorithmic}
	  \REQUIRE $D \in \mathrm{M}_{n\times m}(\mathbb{R})$, $I \in \mathbb{N}^m$, $O \in \mathbb{N}^n$, $\beta \in \mathbb{R}_+$
		\ENSURE $S \in \mathrm{M}_{n\times m}(\mathbb{N})$
		\STATE $S_{ij} \leftarrow 0$
		\WHILE {$\sum_{i=1}^n O_i > 0$}
			\STATE Simulate $i \sim \mathcal{U}_A$ where $A=\left\{k | k\in |[1,n]|,\,O_k \neq 0\right\}$
			\STATE Simulate $j$ from  $|[1,m]|$ with a probability:
					$$P_{i\rightarrow j}=\frac{I_jf(d_{ij},\beta)}{\sum_{k=1}^m I_kf(d_{ik},\beta)}$$
			\STATE $S_{ij} \leftarrow S_{ij}+1$		
			\STATE $I_j \leftarrow I_j-1$ 
			\STATE $O_i \leftarrow O_i-1$	
		\ENDWHILE
	\end{algorithmic}
	{\vspace*{-0.2cm}\hrulefill}	
\end{algorithm}

In \cite{Gargiulo2012}, the authors use a deterrence function $f(d_{ij},\beta)$ with a power law shape: 
\begin{equation}\label{PW}
f(d_{ij},\beta)= d_{ij}^{-\beta} \quad 1\leq i,j \leq n  \enspace .
\end{equation}

\section*{Statistical tools}

This section presents the tools used to calibrate the model and to compare various implementation choices.

\subsection{Calibration of the $\beta$ value.}
\label{Calib}

The same method used in \cite{Gargiulo2012} is used to calibrate the $\beta$ value. $\beta$ is calibrated so as to minimize the average Kolmogorov-Smirnov distance between the simulated commuting distance distribution and one building from the observed data.  
For the basic model we compute the commuting distance distribution with the commuting distance of individuals who are commuting from the region to the region. For the model focused on the outside we compute the commuting distance distribution with the commuting distance of the individuals who are commuting from the region to the region and outside.      

Since the model is stochastic, the final calibration value we consider is the average $\beta$ value over ten replications of the generation process. 

\subsection{An indicator to assess the change.}

It is necessary to have an indicator to compare the simulated commuting network and the observed commuting network (data set 1 in section \ref{observeddata}). Let $R\in \mathrm{M}_{n_1\times n_2}(\mathbb{N})$ represent the observed commuting network when $R_{ij}$ represents the number of commuters from municipality $i$ to municipality $j$. Let $S\in \mathrm{M}_{n_1\times n_2}(\mathbb{N})$ represent a simulated commuting network for the same municipalities. We can calculate the number of common commuters between $R$ and $S$ (Eq. \ref{NCC}) and the number of commuters in $R$ (Eq. \ref{NC}):

\begin{equation}
NCC_{n_1\times n_2}(S,R)=\sum_{i=1}^{n_1}\sum_{j=1}^{n_2}\min(S_{ij},R_{ij}) \enspace
\label{NCC}
\end{equation}
\begin{equation}
NC_{n_1\times n_2}(R)=\sum_{i=1}^{n_1}\sum_{j=1}^{n_2} R_{ij} \enspace
\label{NC}
\end{equation}

From (Eq. \ref{NCC}) and (Eq. \ref{NC}) we calculate the S{\o}rensen similarity index \cite{Sorensen1948}. This index is suitable because it corresponds to the common part of commuters between $R$ and $S$. Thus it is called the common part of commuters (CPC) (Eq. \ref{CPC}):

\begin{equation}
CPC_{n_1\times n_2}(S,R)=\frac{2NCC_{n_1\times n_2}(S,R)}{NC_{n_1\times n_2}(R)+NC_{n_1\times n_2}(S)} \enspace
\label{CPC}
\end{equation}

This index has been chosen for its intuitive explanatory power, as it is a similarity coefficient that provides the likeness degree between two networks. The index ranges from a value of zero, for which there are no any commuter flows in common in the two networks, to a value of one, when all commuter flows are identical between the two networks.

\section{Generating commuting networks for French regions at municipality level}

\subsection{How to cope with regions that are not islands or those that lack detailed data?} 

A commuting network is defined by an origin-destination table (light gray table in Table \ref{ODTable}). At the regional level, this means that it is necessary to know, for each municipality of residence and for each municipality of employment, the value for the flow of commuters traveling from one to another. This kind of data is not always provided by statistical offices and the datasets are usually aggregated: only the total number of out-commuters and in-commuters for each municipality is available for each (dark gray row and colum in Table \ref{ODTable}). To apply the model and define the commuting network, unless we are on a significantly isolated region \footnote[2]{an island for example, in this case gray rows and colums in Table \ref{ODTable} would not exist}, we need to find a way to isolate from the total number of in(out)-commuters (dark gray row and colum in Table \ref{ODTable}) the fraction that relates strictly to the region (light gray table in Table \ref{ODTable}). However, this is not a simple task.

Furthermore, even if these parts can be isolated, a problem remains due to the border effect. Indeed, if we consider only the region, there is the risk of making an error in the reconstruction of the network for municipalities near the region's border. The higher the proportion of individuals working outside of the region, the more significant the error will be.

To go further, we propose to change the inputs for the algorithm. Instead of only considering the regional municipalities as possible places of work, we also consider an $outside$ of the region. The outside represents the surroundings of the studied area. The following section describes a method for considering this outside area practically.

\begin{table*}
	\caption[Origin-destination table]{Origin-destination table; The light gray table represents the commuters living and working in the region for each municipality of the region; The gray column represents the out-commuters living in the region and working outside (Out.) for each municipality of the region; The gray line represents the in-commuters working in the region and living outside (Out.) for each municipality of the region; The dark gray line(column) represents the total number of out(in)-commuters for each municipality of the region.}
	\label{ODTable}
		\begin{center}
    \begin{tabular}{c}
			\includegraphics[scale=0.35]{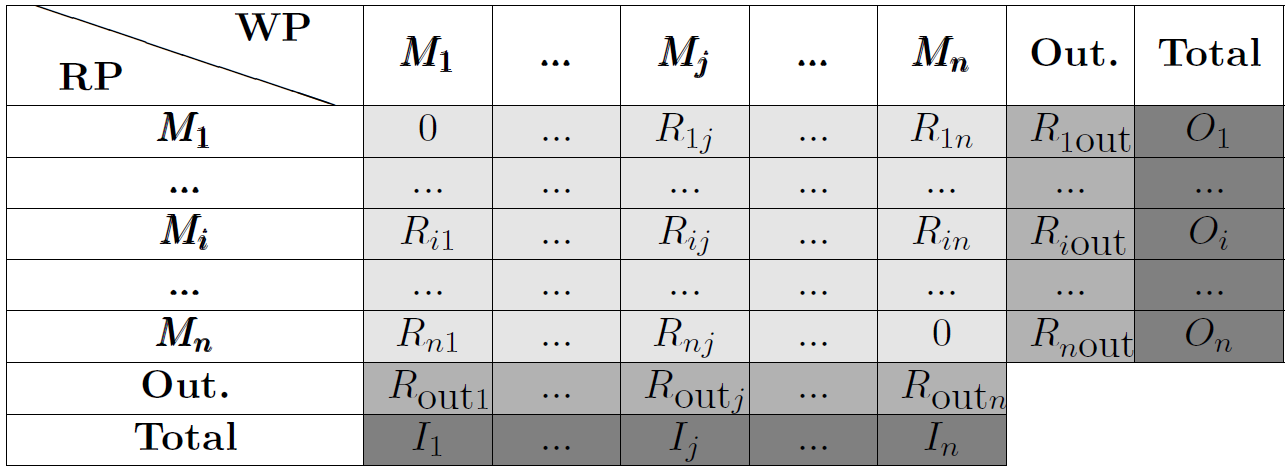}
		\end{tabular}
	\end{center}
\end{table*}

\begin{table*}
	\caption[Origin-destination table from the region to the region and the outside]{Origin-destination table from the region to the region and the outside; The light gray table represents the commuters living (place of residence RP) and working (place of work WP) in the region for each municipality of the region; The gray table represents the commuters living (place of residence RP) in the region and working (place of work WP) outside of the region.}
\label{netout}
		\begin{center}
    \begin{tabular}{c}
			\includegraphics[scale=0.35]{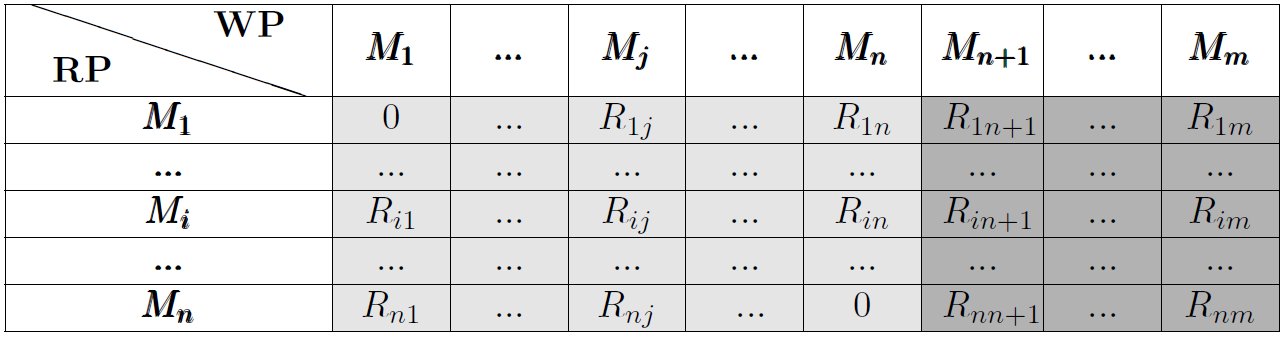}
		\end{tabular}
	\end{center}
\end{table*}

\subsubsection{A new extended to outside job search base. }

We implement the model to generate 23 various regions in France. Their outside is composed of the set of municipalities of their neighboring "departments".  

We consider the outside of the region to be composed of $m-n$ municipalities, where $n$ represents the number of municipalities in the region. The inputs are the directly available aggregated data at the municipal level: 

\begin{itemize}
\item $D=(d_{ij})_{1 \leq i \leq n \atop 1 \leq j \leq m}$ the Euclidean distance matrix between municipalities both in the same region and in the outside.
\item $(I_j)_{1 \leq j \leq m}$ the total number of in-commuters of municipality $j$ of the region and outside of it (i.e. the number of individuals working in municipality $j$ of the region or the outside and living in another municipality).
\item $(O_i)_{1 \leq i \leq n}$ the total number of out-commuters of municipality $i$ of the region only (i.e. the number of individuals living in  municipality $i$ of the region and working in an other municipality).
\end{itemize}

The purpose of the algorithm that introduces the outside is to generate the origin-destination table (light gray and gray sub-table in Table \ref{ODTable}). To do this the algorithm presented in Algorithm \ref{commut} is used to simulate the Table \ref{netout}. From this, through difference the Table \ref{ODTable} can be obtained with the total number of in-commuters $(I_j)_{1 \leq j \leq n}$, the total number of out-commuters $(O_i)_{1 \leq i \leq n}$ and the light gray table of the Table \ref{netout}.

A matricial representation of the origin-destination table presented in the light gray and gray sub-table in Table \ref{ODTable}, known as the simulated matrix $S \in \mathrm{M}_{(n+1)\times (n+1)}(\mathbb{N})$ is obtained. $S_{ij}$ represents: 
\begin{itemize}
\item the number of commuters from municipality $i$ (in the region) to municipality $j$ (in the region) if $i,j \neq n+1$; 
\item the number of commuters from outside to municipality $j$ (in the region) if $i=n+1$ and $j\neq n+1$; 
\item the number of commuters from municipality $i$ to outside if $i \neq n+1$ and $j=n+1$. 
\end{itemize}

\subsubsection{Comparison of the two models: Assessing the impact of the outside. }

We assess the impact of the outside through a comparison between the network generations for 23 French regions both with and without the outside. The generation is made on a municipality scale using a power law deterrence function. 

\begin{figure}
   \includegraphics[width=\linewidth]{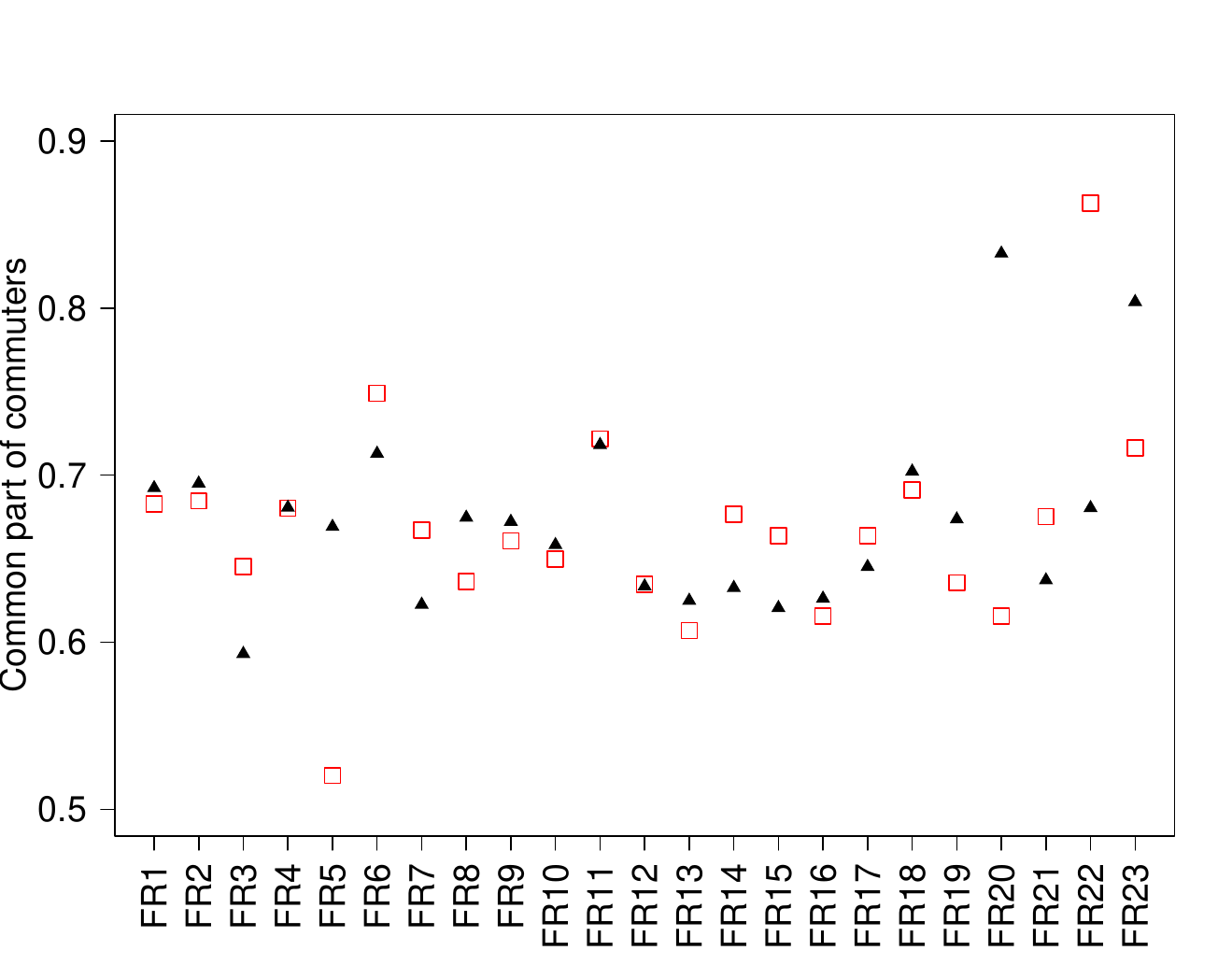}
    \caption{Average CPC for 23 regions. The squares represent the basic model; The triangles represent the model with outside.}
    \label{Fig1}
\end{figure}

Both implementations are compared through their CPC values between the simulated network $S$ and the observed network $R$ (data set 1 presented in Section \ref{observeddata}) for each region. We replicate the generation for each region ten times and our indicator on each replicate is calculated. In all the presented figures, the indicator averages ten replications. The variation of the indicator over the replications is very low, averaging $1.02\%$ at most. Consequently, this is not represented on the figures. Fig. \ref{Fig1} presents the common part of commuters $CPC_{n\times n}(S,R)$ between the simulated network $S$ and the observed network $R$; The squares represent the CPC between the observed network $R$ and the simulated networks obtained with the regional job search base; The triangles represent the CPC between the observed network $R$ and the simulated networks obtained with a job search base comprising the region and its outside. It's important to note that for the implementation without outside $S\in \mathrm{M}_{n\times n}(\mathbb{N})$ while for the implementation with outside $S\in \mathrm{M}_{(n+1)\times (n+1)}(\mathbb{N})$. In order to compare the two models, the regional network (commuters from the region to the region) must be taken into consideration. Indeed, in the without-outside cases $NC_{n\times n}(S)=NC_{n\times n}(R)$ but this is not necessarily true for the with-outside cases.

Fig. \ref{Fig1} shows that the two job search bases give results which are not different. Thus, introducing the outside solves the problem linked to a lack of detailed data without changing the quality of the resulted simulated network. Indeed, one must keep in mind that the inputs for the with-outside cases do not require detailed data in comparison to the without-outside cases. 

\subsection{Choosing a shape for the deterrence function}

The next problem relates to the form of the deterrence function which rules the impact of distance on the choice of the place of work relative to the quantity of job offers. The initial work done by \cite{Gargiulo2012} proposes to use a power law. However, \cite{Barthelemy2011} states the form of the deterrence function varies significantly, and can sometimes be inspired by an exponential function as in  \cite{Balcan2009} or by a power law function as in \cite{Viboud2006}. Through choosing the much more convenient deterrence function, we compare the quality of generated networks for 34 French regions obtained with the model with outside using both the exponential law and the power law. 

A deterrence function following an exponential law is introduced: 
\begin{equation}\label{exp}
    f(d_{ij},\beta)= \displaystyle e^{-\beta d_{ij}} \quad 1\leq i \leq n \mbox{ and } 1\leq j \leq m \enspace .
\end{equation}

To compare the two deterrence functions, we have generated the networks of 34 various French regions (see Table \ref{tabdata} for details) that replicate ten times for each region. The networks were generated with a job search base for the algorithm that considers the outside. 

\begin{figure}
    \includegraphics[width=\linewidth]{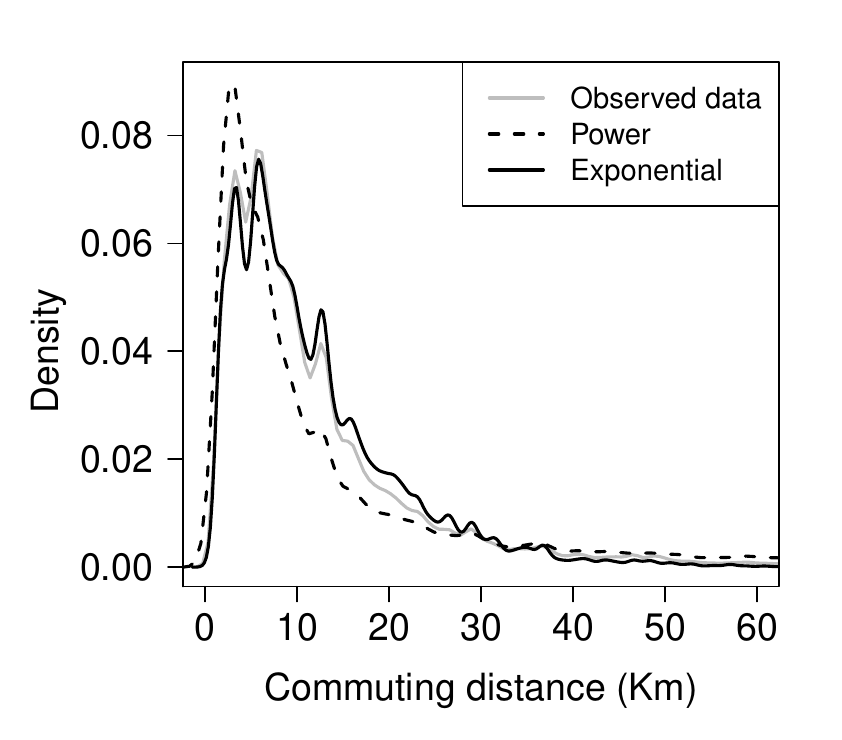}
    \caption{Density of the Auvergne commuting distance distribution; the solid line represents the observed commuting distance distribution; the dotted line represents the commuting distance distribution obtained with the calibrated model with a job search base comprising the outside and the exponential law; the dashed line represents the commuting distance distribution obtained with a job search base comprising the outside and the power law. The two simulated commuting distance distributions are computed for one replication each.}
   \label{Fig2}
\end{figure}

For example, Fig. \ref{Fig2} shows that we obtained a better estimation of the Auvergne commuting distance distribution when using the exponential law. We computed the observed commuting distance distribution with the observed Auvergne commuting network (data set 1 presented in Section \ref{observeddata}) and the Euclidean distances between the Auvergne municipalities (data set 3 presented in Section \ref{observeddata}).  

\begin{figure}
   \includegraphics[width=\linewidth]{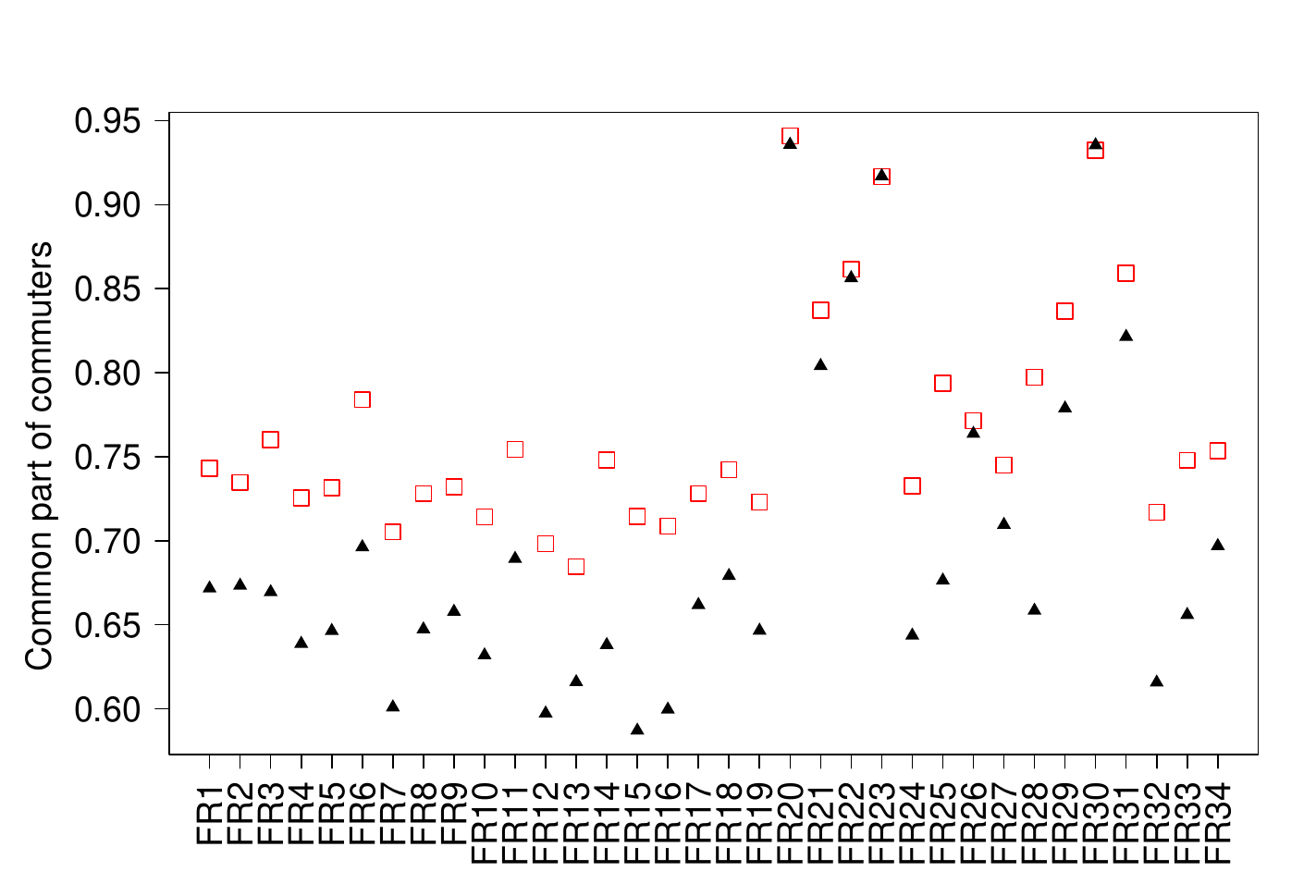} 
   \caption{Average CPC for the power shape (triangle) and the exponential shape (square) for 34 French regions.}
   \label{Fig3}
\end{figure}

More systematically, we plot, for the exponential law and power law, the average of the replications for the common part of commuters $CPC_{(n+1)\times (n+1)}(S,R)$ between the simulated network $S$ and the observed network $R$ in Fig. \ref{Fig3}. This clearly indicates that the average proportion of common commuters is always better when using an exponential law represented by squares. 

\begin{figure*}
  \centering  
  \includegraphics[scale=0.21]{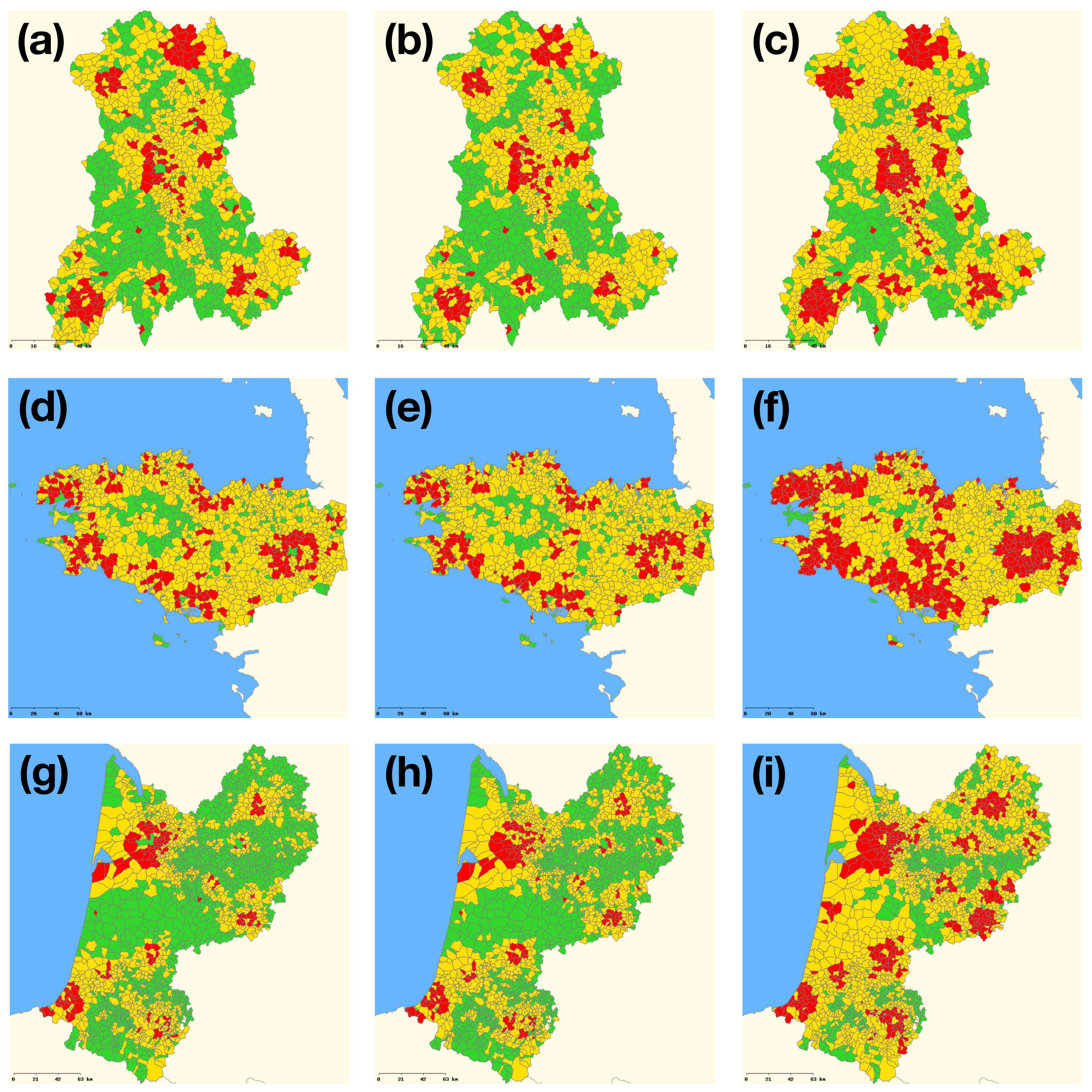}
  \caption[Maps of the average CPC by municipalities]{Maps of the average CPC by municipalities obtained with ten replications. In green CPC $\leq0.5$; In yellow $0.5<$ CPC $\leq0.75$; In red $0.75<$ CPC. (a), (d) and (g) Model with the power shape without outside; (b),(e) and (h) Model with the power shape with outside; (c), (f) and (i) Model with the exponential shape with outside. (a)-(c) Auvergne case-study; (d)-(f) Bretagne case-study; (g)-(i) Auquitaine case-study. \textit{\scriptsize{Base maps source: Cemagref - DTM - 
  D\'{e}veloppement Informatique Syst\`{e}me d'Information et Base de Donn\'{e}es : F.Bray \& A.Torre
IGN (G\'{e}ofla\textsuperscript{{\fontsize{5}{5}\textregistered}}, 2007).}}}
  \label{Fig5}
\end{figure*}

\subsection{Spatial Analysis}

To better understand how CPC is spatially distributed at a more granular level we mapped the CPC by municipality for three models and three study areas. In Fig. \ref{Fig5}, it can be observed that for all case studies (in rows) the highest values of the CPC were obtained by municipalities using the model with an exponential shape including the outside (third column). It can also be noted that the model without the outside (second column) and the model with the power shape including the outside (first column) give results which are not wholly different.

As we can see in Fig. \ref{Fig5}, the CPC values are not uniformly distributed in the municipalities of the three areas. The error seems to increase as distance from the urban areas increases.

\begin{figure*}
	\includegraphics[width=\linewidth]{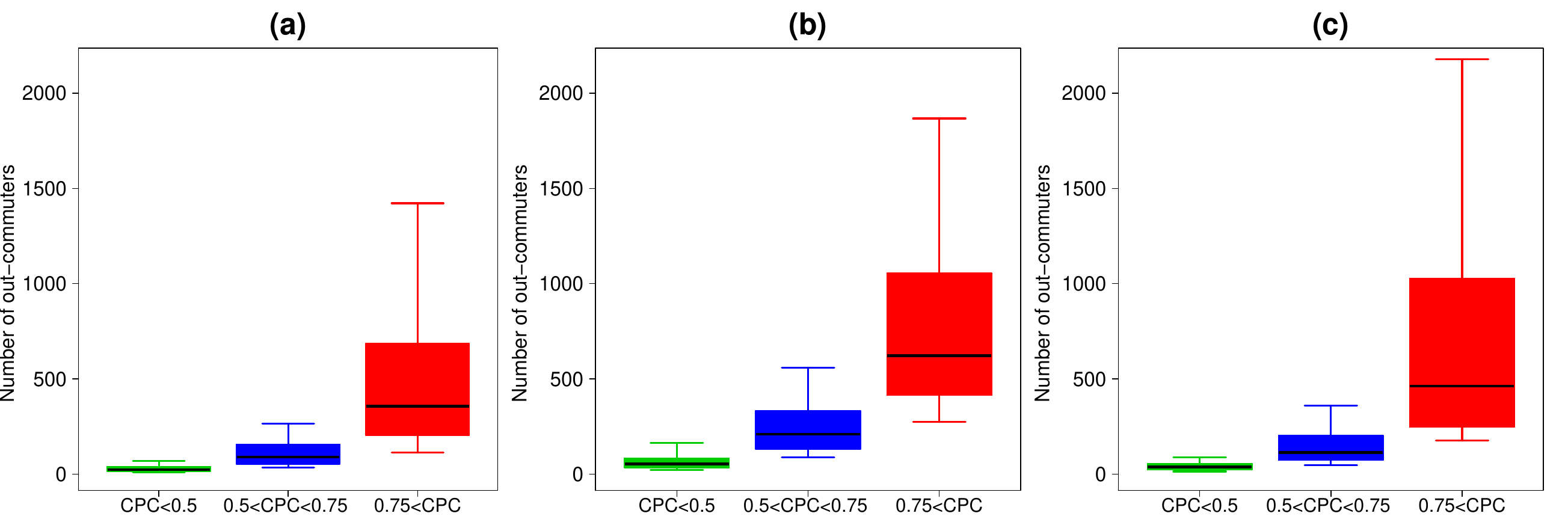}
	\caption{Boxplots of the number of out-commuters in term of the CPC by municipality for the model with the exponential shape with outside. (a) Auvergne case study; (b) Bretagne case study; (c) Aquitaine case study.}
	\label{Fig4}
\end{figure*}

We now focus on the third model with an exponential shape including the outside to better understand which types of municipalities compose the three clusters (CPC$\leq0.5$, $0.5<$CPC$\leq0.75$ and $0.75<$CPC). We identify the number of out-commuters as the most explanatory variable. Indeed, we can observed in Fig. \ref{Fig4} that the distribution of the number of out-commuters in each cluster is significantly different. The higher the average number of out-commuters, the higher the CPC. Having performed analyses of variance (ANOVA) for each case study, we obtained significant differences between the averages for the number of out-commuters in each cluster with a 0.95\% level of confidence for each case study.  

For the three regions, the CPC value is strongly linked to municipality characteristics. Indeed, the municipalities with $0.75<$CPC are urban and suburban municipalities with a high number of out-commuters that are closed to a large urban municipality. In contrast, the municipalities with a low number of out-commuters that are far from large urban municipalities have a CPC lower than 0.5. For this type of municipality, the commuting flows are very small. Thus they are  difficult to reproduce with the mechanisms taken into consideration. However, the distance to cities does not appear to be particularly responsible for the error. The timing for the job offer arrival on the job market is probably much more significant in determining the local topology of the network than elsewhere. These flows represent about 4\% of the total number of out-commuters for the Auvergne region, 1\% for Bretagne and 5\% for Aquitaine.

\subsection{Calibrating the model for French regions}

The final problem involves the calibration process, which previously required detailed and accurate data.\\

\begin{figure}
			\includegraphics[width=\linewidth]{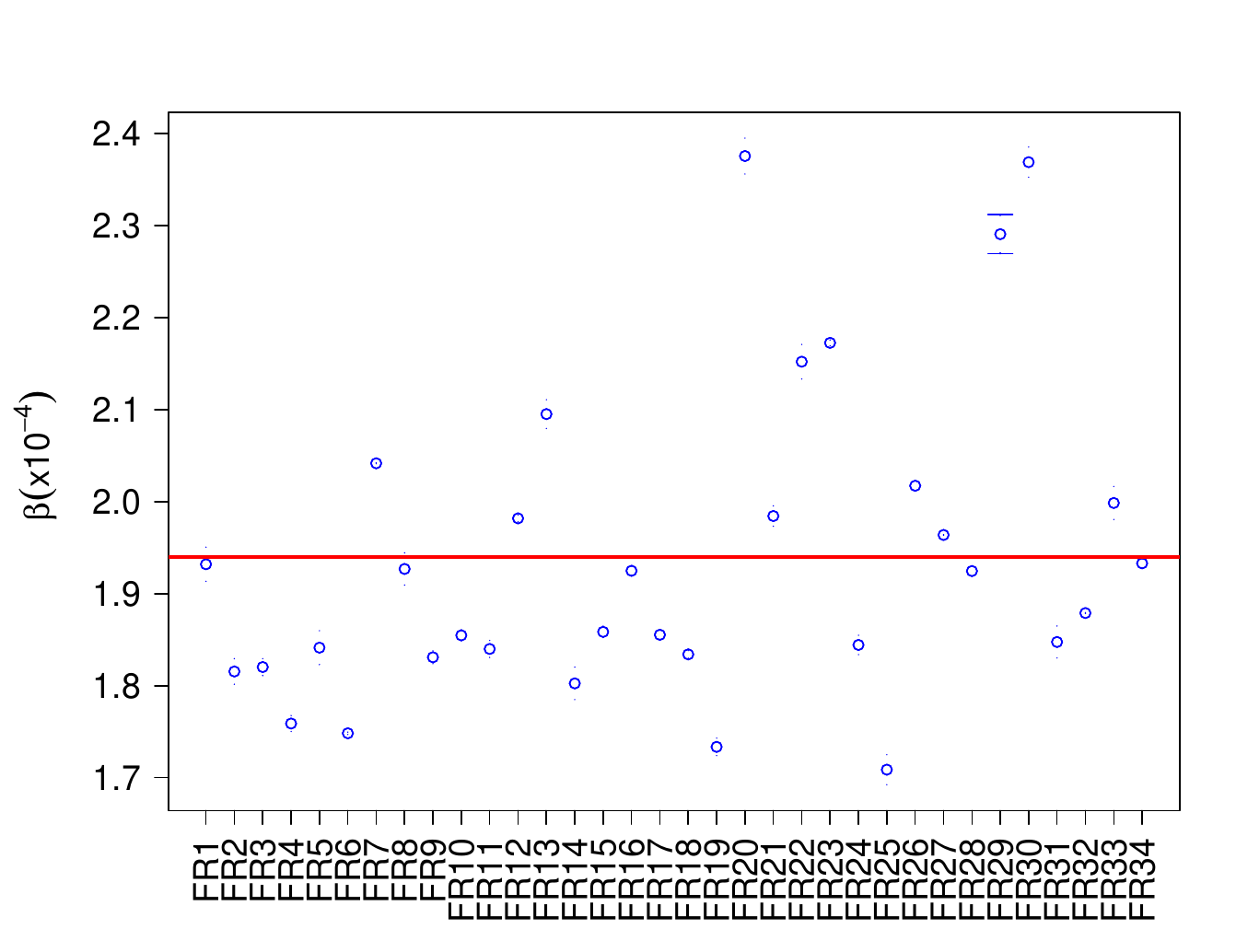}
			\caption{The circle represents the average calibrated $\beta$ values for ten replications (The confident interval is composed of the minimum and 
			the maximum) for each regions; the line represents the average $\beta$ value for the $34$ regions.}
		  \label{Fig6} 
\end{figure}

Fig. \ref{Fig6} shows the calibrated $\beta$ values for each of the $34$ regions in France. It can be observed that these values display subtle variations from about $1.7\cdot10^{-4}$ to $2.4\cdot10^{-4}$ with the average $\beta$ valued ($C=1.94\cdot10^{-4}$) corresponding to the dark line. 

\begin{figure}
	\centering
	\includegraphics[width=\linewidth]{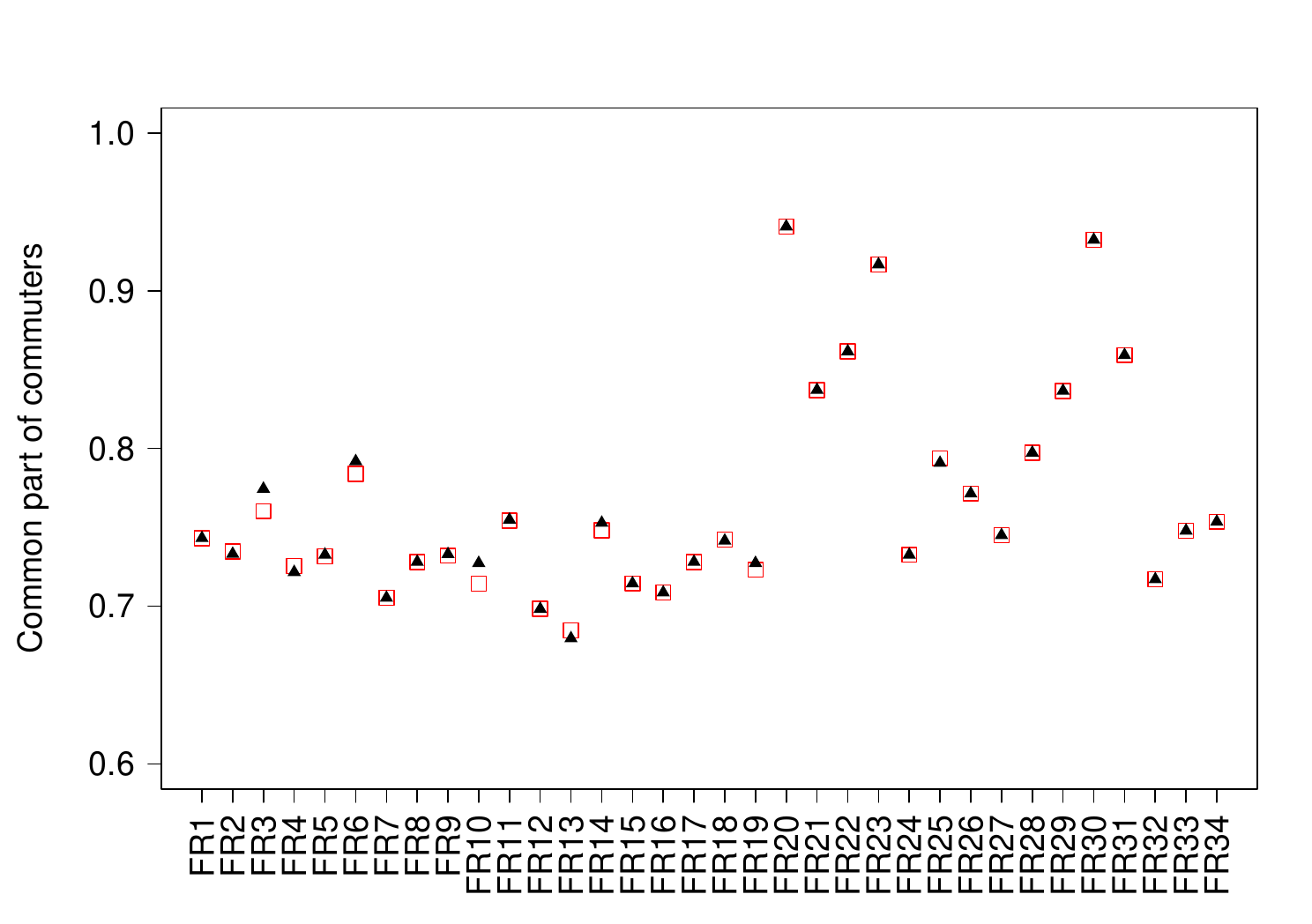}
	\caption{Common part of commuters for the $34$ regions; The squares represents the average CPC (10 replications) obtained with the calibrated $\beta$ value; The triangles represents the average CPC (10 replications) obtained with the estimated $\beta$ values (average $\beta$ value over the 34 calibrated $\beta$ values).}
	\label{Fig7}
\end{figure}

Then we hypothesize that it is possible to directly calibrate the algorithm to generate the $34$ regions in France, by using a constant equal to $C$. To study the influence of this approximation on the common part of commuters we have computed the CPC with $C$ as the parameter value for the $34$ regions. We observe in Fig. \ref{Fig7} that the influence of the $\beta$'s approximation on the CPC is very weak. It can then be noted that the average CPC obtained with $C$ is, for some regions, higher than the CPC obtained by the $\beta$ value that is not averaged. It is possible that the common part of commuters is better with another beta value because it is not a calibration criterion.  

It is not necessary to study the influence of the $\beta$'s approximation on the calibration criterion. Indeed, from the studies made by \cite{Gargiulo2012}, we know the CPC and the calibration criterion show a significant correlation. The CPC and the calibration criterion follow the same evolution in terms of $\beta$. The $\beta$ value for minimization of the Kolmogorov-Smirnov distance is very close to the one obtained for maximization of the CPC (see the figure 7 in \cite{Gargiulo2012} which perfectly illustrates this relation). The CPC values remain quasi-identical to $\beta$=$C$ or to $\beta$ valued from the calibration process presented in Section \ref{Calib}, the quality of the approximation of the calibration criterion, i.e. the commuting distance distribution, remains the same.

\section{Discussion}

To study the rural area dynamics through micro-simulation, we need virtual commuting networks that link individuals living in the municipalities of various French regions. As the studied scale is very low, the flows are low, and we thus decided to opt for a stochastic generation algorithm. The one recently proposed by \cite{Gargiulo2012} is relevant to our problem. Starting from this model, we implement the commuting networks of 34 different French regions. The implementation work leads us to solve three practical problems. 

The first problem involves the fact that our French regions are not islands. Indeed, some of the inhabitants, especially those living close to the border of the region, are likely to work in municipalities located outside the region of residence. However, classical approaches to generating commuting networks consider only residents of the region that work in the region. That is also the case for ours. Data providing details, or knowledge, allowing the modeler to evaluate people living in the region but working outside is difficult to obtain. Thus, we address this issue by extending the geographical base of the job search for commuters living in the municipalities to a sufficiently large number of municipalities located outside the region of residence. We compare the model without municipalities located outside and the model with outside municipalities to 23 French regions. We are able to come to a conclusion regarding the relevance of our solution which keeps the value of our quality indicator identical. At the same time, it is not necessary to have information regarding those who do not work in the region, which allows us to generate networks using only the aggregated data.

The \cite{Gargiulo2012} model is based on the gravity law. Then, our second problem relates to the deterrence function, which is more of a power law or an exponential law depending on the study. Moreover, as empirical studies comparing generated networks to observed data are extremely rare \cite{Barthelemy2011}, few know which is better. In order to select the more convenient one for our French regions, we have compared the quality of generated networks for 34 regions obtained with both the exponential law and the power law.  Better results were obtained with the exponential law, no matter the region. Indeed, the 34 regions display significant variance in regards to surface area, the number of municipalities, and the number of commuters.

The final problem involved calibration. Applying a model with an extended job search base and an exponential deterrence function, we found a constant equal to $1.94\cdot10^{-4}$ to be a perfect parameter value for generating commuting networks for French administrative regions, no matter the region. However, we did not test this result for other countries with different types of administrative regions. The robustness of this result to commuting networks of different scales has been studied in \cite{Lenormand2012b}. The $\beta$ value correlated to a scale consistent with the results obtained in this paper.

A spatial analysis of three different case studies has been proposed, and it was shown that the CPC value by municipality strongly correlated with the number of out-commuters for the municipality. Our model is not able to reproduce very small flows which represent between 1 and 5\% of the total flows in the region we studied. However, we continue to question if it makes sense to attempt to reproduce them.

\section*{Acknowledgements}

This publication has been funded by the Prototypical policy impacts on multifunctional activities in rural municipalities collaborative project, European Union 7th Framework Programme (ENV 2007-1), contract no. 212345. The work of the first author has been funded by the Auvergne region.

\bibliographystyle{unsrt}  
\bibliography{JTLU}

\begin{table*}
\caption{Description of the regions}
\label{tabdata}
\begin{center}
\begin{tabular}{>{\centering}m{1cm}>{\centering}m{3cm}>{\centering}m{2cm}>{\centering}m{2cm}>{\centering}m{2cm}>{\centering}m{2cm} m{2cm}<{\centering}}
\hline
\textbf{ID} & \textbf{Region} & \textbf{Number of municip. (region)} &  \textbf{Number of municip.  (outside)} & \textbf{Region area (km$^2$)} & \textbf{Average municip. area (km$^2$)}  & \textbf{Number of commuters}\\
\hline
FR1	&	Auvergne	&	1310	&	3463	&	26013	&	19.86	&	295776	\\
FR2	&	Bretagne	&	1269	&	1447	&	27208	&	21.44	&	653710	\\
FR3	&	Ain	&	419	&	2809	&	5762	&	13.75	&	162370	\\
FR4	&	Alsace	&	903	&	3081	&	8280	&	9.17	&	440961	\\
FR5	&	Aquitaine	&	2296	&	2835	&	41309	&	17.99	&	700452	\\
FR6	&	Mayenne	&	261	&	3124	&	5175	&	19.83	&	69915	\\
FR7	&	Loz\`ere	&	185	&	1859	&	5167	&	27.93	&	12273	\\
FR8	&	Poitou-Charente	&	1464	&	2467	&	25810	&	17.63	&	375363	\\
FR9	&	Centre	&	1842	&	4718	&	39151	&	21.25	&	624693	\\
FR10	&	Midi-Pyr\'en\'ee	&	3020	&	3845	&	45348	&	15.02	&	546162	\\
FR11	&	Limousin	&	747	&	3169	&	16942	&	22.68	&	139481	\\
FR12	&	Franche-Comt\'e	&	1786	&	3317	&	16202	&	9.07	&	268399	\\
FR13	&	Haute-Normandie	&	1420	&	3536	&	12317	&	8.67	&	469335	\\
FR14	&	Haute-Marne	&	433	&	3914	&	6211	&	14.34	&	42690	\\
FR15	&	Vosges	&	515	&	3808	&	5874	&	11.41	&	92053	\\
FR16	&	Lorraine	&	2339	&	3067	&	23547	&	10.07	&	547457	\\
FR17	&	Creuse	&	260	&	1814	&	5565	&	21.40	&	23949	\\
FR18	&	Languedoc-Roussillon	&	1545	&	3046	&	27367	&	17.71	&	409116	\\
FR19	&	Charente-Maritime	&	1948	&	1983	&	25606	&	13.14	&	375363	\\
FR20	&	Haut-de-Seine	&	36	&	1245	&	176	&	4.89	&	973173	\\
FR21	&	Yveline	&	262	&	1543	&	2284	&	8.72	&	618741	\\
FR22	&	Val d'Oise	&	185	&	1707	&	1246	&	6.74	&	526600	\\
FR23	&	Val de Marne	&	47	&	1234	&	245	&	5.21	&	642092	\\
FR24	&	Haut-Rhin	&	377	&	2283	&	3525	&	9.35	&	183504	\\
FR25	&	Tarn et Garonne	&	195	&	2338	&	3718	&	19.07	&	41600	\\
FR26	&	Pyr\'en\'ee-Atlantique	&	547	&	449	&	4116	&	7.52	&	65469	\\
FR27	&	Alpes-Maritimes	&	163	&	353	&	4299	&	26.37	&	163445	\\
FR28	&	Loire	&	327	&	2788	&	4781	&	14.62	&	178828	\\
FR29	&	Territoire de Belfort	&	102	&	2031	&	609	&	5.97	&	45185	\\
FR30	&	Seine-Saint-Denis	&	40	&	783	&	236	&	5.90	&	655200	\\
FR31	&	Essonne	&	196	&	1597	&	1804	&	9.20	&	518321	\\
FR32	&	Ardennes	&	463	&	2588	&	5229	&	11.29	&	59963	\\
FR33	&	Aube	&	433	&	2728	&	6004	&	13.87	&	75561	\\
FR34	&	Corr\'eze	&	286	&	2088	&	5857	&	20.48	&	49815	\\
\hline   
\end{tabular}
\end{center}
\end{table*}

\end{document}